\documentclass[preprint,12pt,authoryear]{elsarticle}
\usepackage{epsfig}

\usepackage{amsmath}
\usepackage{amssymb}

\newcommand{\bx}{{\bf x}}
\newcommand{\coeff}{a}
\newcommand{\bcoeff}{{\bf a}}
\newcommand{\mathel}{{\mathcal L}}

\newcommand{\bara}{{\bar {A}}}
\newcommand{\barb}{{\bar {B}}}
\newcommand{\barc}{{\bar {C}}}
\newcommand{\bard}{{\bar {D}}}
\newcommand{\fft}{{\mathcal F}}

\newcommand{\mathk}{{\mathcal K}}

\newcommand{\bu}{{\bf u}}
\newcommand{\eye}{\,{\bf I}}

\newcommand{\ft}{\,{{\mathcal F}}}
\newcommand{\bnabla}{{\boldsymbol \nabla}}

\journal{Journal of Computational Physics}

\begin{document}
\begin{frontmatter}

\title{Spatio-spectral concentration of Convolutions}


\author{Shravan M. Hanasoge}
\address{Department of Astronomy and Astrophysics, Tata Institute of Fundamental Research, Mumbai 400005, India}
\ead{hanasoge@tifr.res.in}




\begin{abstract}
Differential equations may possess coefficients that vary on a spectrum of scales. 
Because coefficients are typically multiplicative in real space, they turn into convolution operators in spectral space, mixing all wavenumbers. 
However, in many applications, only the largest scales of the solution are of interest and so the question turns to whether
it is possible to build effective coarse-scale models of the coefficients in such a manner that the large scales of
the solution are left intact. 
Here we apply the method of numerical homogenization to deterministic linear equations to generate sub-grid-scale models
of coefficients at desired frequency cutoffs.
We use the Fourier basis to project, filter and compute correctors for the coefficients.
The method is tested in 1D and 2D scenarios and found to reproduce the coarse scales of the solution to varying
degrees of accuracy depending on the cutoff. We relate this method to mode-elimination Renormalization Group (RG) and discuss the connection between accuracy and the cutoff wavenumber. The tradeoff is governed by a form of the uncertainty principle for convolutions, which states that as the convolution operator is squeezed in the spectral domain, it broadens in real space. As a consequence, basis sparsity is a high virtue and the choice of the basis can be critical. 
\end{abstract}


\end{frontmatter}


\section{Introduction}
\label{intro}
In the early 1970s, a significant hurdle faced by particle physicsists was the computation of partition functions that involved evaluating integrals over large ranges in wavenumbers. Ultimately, the only parameters of interest were at coarse scales of the systems under study. 
To compute these parameters, Kenneth Wilson introduced mode-elimination Renormalization Group (RG) that enabled building low-wavenumber representations of fluctuating coefficients while succeeding in preserving the coarse-scale accuracy of the solutions. In other words, RG describes a means of projecting the small scales onto the large, and by degrees, integrating out rapid variations in the coefficients. This procedure has the potential to greatly reduce the computational burden. RG is perhaps the most celebrated instance of the concept of building coarse-scale models of inherently multi-scale phenomena, and has been widely used to model a range of phenomena. 
The application of mode-elimination RG to generate sub-grid-scale models of fluid turbulence, i.e. that convey the effect of scales smaller than the grid size, was suggested by \cite{yakhot86} and studied in detail for passive-scalar advection by \cite{majda90} \citep[also see \cite{smith98} who compare these methods and Kraichnan's Direct Interaction Approximation,][]{kraichnan62_DIA}. \cite{yakhot86} suggested integrating over shells of wavenumbers, proceeding from the largest to the smallest, sequentially adding corrections to the coefficients of the turbulence model in question.

Independently, \cite{kozlov79} and \cite{varadhan82} studied solutions to the diffusion equation with random coefficients. Defining a small parameter $\varepsilon$ to be the ratio of the correlation length-scale of the random coefficent to the relevant coarse scale of the solution, \cite{kozlov79} and \cite{varadhan82} derived a two-scale asymptotic theory to estimate effective coefficients. The formulation is similar in concept to RG and produces a zero-wavenumber (constant) representation of the original random coefficient, and in the limit of vanishingly small $\varepsilon$, \cite{varadhan82} showed that the solution was accurate. Multi-scale coefficients however present a challenge and appropriately decimating these coefficients onto a coarser grid (not necessarily zero wavenumber) is of relevance in computational physics. In the finite-scale-cutoff scenario with non-random media, the method of numerical homogenisation is remarkably similar to RG. These techniques can be applied to build effective models of wavespeeds (whose spatial distributions can be complicated) in seismology, deriving coarse-scale descriptions of porosity coefficients in porous media etc. In this article, we will focus on these two scenarios.

A central goal of seismic studies of the Sun, stars and Earth is to infer the structural and dynamical properties of interiors using observations of their surface oscillations. The {\it forward problem}, critical to this effort, is the simulation of small-amplitude (linear) waves through the relevant media. Such media can comprise a wide spectrum of length scales, possibly much smaller than the wavelength.
The problem thus becomes computationally stiff and very expensive to attempt. We are therefore interested in bringing to bear methods of homogenization, which describe the coarse-scale behavior of differential equations with rapidly varying coefficients, on these problems of wave propagation. We seek to replace the fine-scale structure with an effective sub-grid-scale model such that the coarse scales of the solution are accurately reproduced to within a specified tolerance.

Similarly, in porous media, the permeability of the medium, a tensor quantity, is finely sampled at a large number of spatial points. The goal is to coarsen the grid and appropriately average these tensor coefficients. Effective coarse-grained models of fine-scale tensor coefficients will necessarily mix various components. Here we will derive a formal theory that describes how to mix various terms.
Classical homogenization primarily addresses problems in which the coefficients periodically vary \citep{lions78}, with the sub-grid model being a zero-wavenumber representation (e.g., by the harmonic mean). However, in a number of real-world applications, the rapid variations are aperiodic and a more general theory to attempt such problems is required. Along these lines, multi-resolution analysis
in the aid of numerical homogenization of aperiodic media has been developed \citep[e.g.][]{brewster, dorobantu, engquist02,owhadi}. More recently, 
e.g., \cite{capdeville1d}, have posed the problem of terrestrial seismic wave propagation through aperiodic heterogeneous media 
in the language of classical homogenization. 

In this article, we follow the methodology of \cite{dorobantu} and \cite{engquist02}, which is well laid out and from which details may be intuited. Numerical homogenization affords two major advantages: significant reduction in spatial complexity and a less restrictive Courant condition on the timestep. Here, we use the spatial Fourier and Haar-wavelet bases to investigate the accuracy of numerical homogenization on three different wave equations, each gaining complexity over the previous. The Fourier basis lends itself to elegant interpretation but produces dense matrices whose inverses may not be easy to compute. In contrast sparse matrix inversion techniques may be easily extended to homogenization in the Haar basis. The demonstrable success of the method encourages a more complete exploration of its possibilities.

\section{Numerical homogenization in 1D}\label{num1d}
Consider the 1D operator $\mathel$ acting on a function $u$ defined by $\mathel u = \partial_x (\coeff\, \partial_x u)$, where $\coeff  = \coeff(x) > 0$ is a coefficient, $\partial_x$ is the spatial derivative with respect to the $x$ coordinate. 
The wave equation corresponds to $\partial^2_t u- S = \mathel u $, where $t$ is time, $S = S(x,t)$ is a source and the equation takes on a hyperbolic character. The time-independent porous-flow equation, identical to the diffusion equation, is given by $\mathel u = 0$ and is elliptic in character.

The product in real space between $\coeff(x)$ and $\partial_x u$ is a convolution in Fourier domain, resulting in the mixing of coarse and fine scales.
In other words, the Fourier transform of this term (we do not add extra symbols to denote the transformed quantity) is 
\begin{equation}
{\mathel u} (k) = - \sum_{k'} k k' \coeff(k-k')\,u(k'),
\end{equation}
and this results in a mixing between low and high wavenumbers. Thus to obtain the low-wavenumber representation of $u$, one must solve the equation over the full set of wavenumbers, which can be computationally expensive.
The goal then is to create a sub-grid-scale model of $\coeff$ such that the coarse scales of $u$ are well reproduced. 
Let us consider the projection of a function in the Fourier basis. Define a projection operator $\ft$ that transforms a function in real space
to the Fourier basis, producing a set of Fourier coefficients which may subsequently be characterized as ``coarse" or ``fine". Denoting the forward transform by $\fft$,
and given an $N\times1$ vector $v$, the projection is written as $\fft v = \begin{pmatrix} P \\Q\end{pmatrix}v$, where $P$ is a $k_p\times N$ matrix that projects
$v$ on to the coarse set of coefficients (of size $k_p \times 1$) and the $(N-k_p)\times N$-sized matrix $Q$ projects $v$ on to the fine coefficients. 
We note that $\fft^{-1}\fft = \fft\fft^{-1}=\eye_N$, where the subscript denotes the size of the identity matrix ($N\times N$).
Since we use the orthogonal Fourier basis, the inverse transform is $\fft^{-1} = (P^*\,\,\,\, Q^*)$, where the $*$ denotes conjugate transpose (the Hermitian transpose) and the associated identities are satisfied 
\begin{equation}
PP^* = \eye_P,~~ QQ^* = \eye_Q,~~PQ^* = {\bf 0},~~QP^* = {\bf 0},~~P^*P + Q^*Q = \eye_N. 
\end{equation}
Note that $\eye_P$ is of size $k_p\times k_p$ and $\eye_Q$ of size $k_q\times k_q$ (where $k_q = N - k_p$).
With no loss of generality, this method may also be extended to other orthogonal and bi-orthogonal systems. Note that the sizes of the zero matrix ${\bf 0}$ in the two identities in the middle are $k_p\times k_q$ and $k_q\times k_p$ respectively. 
We would like to obtain the effective coarse-scale representation of this spatial operator. 
First we consider the projection of the differential operator on the Fourier basis. The Fourier transform of { $\partial_x u(x)$ is $ik\,u(k)$, where $u(k)$ denotes the Fourier transform
of $u(x)$ (note that, as before, we do not use additional symbols to denote transformed quantities).} 
Thus the Fourier projection of the spatial derivative operator is a diagonal matrix which we denote as
\begin{equation}
\begin{pmatrix} P \\ Q\end{pmatrix}\partial_x\begin{pmatrix} P^* & Q^*\end{pmatrix} = \begin{pmatrix} K_P & 0 \\ 0 & K_Q\end{pmatrix},
\end{equation}
where diagonal matrices (referred to by `diag' in the equation below)
\begin{eqnarray}
K_P &=& {\rm diag}(ik) \,\, : \,\,|k| < k_p,\\
K_Q &=& {\rm diag}(ik) \,\, : \,\,|k| \ge k_p,
\end{eqnarray}
contain low ({\it coarse}) and high ({\it fine scale}) wavenumbers respectively. The zero-wavenumber (i.e., the mean or the dc component of $u$) is
an element in $K_P$, which means that $K_P$ is non-invertible whereas $K_Q$ possesses an inverse.
Thus $\partial_x u$ projected on to the Fourier basis is
\begin{equation}
\begin{pmatrix} P \\ Q\end{pmatrix}\partial_x u = \left[\begin{pmatrix} P \\ Q\end{pmatrix}\partial_x\begin{pmatrix} P^* & Q^*\end{pmatrix} \right]\begin{pmatrix} P \\ Q\end{pmatrix} u = 
 \begin{pmatrix} K_P & 0 \\ 0 & K_Q\end{pmatrix} \begin{pmatrix} u_P \\  u_Q\end{pmatrix}.
\end{equation}

Similarly the projection of the full term $\partial_x\,( \coeff\, \partial_x u)$
on to the Fourier basis is given by
\begin{equation}
{\mathel u}(k) =  \begin{pmatrix} P \\ Q\end{pmatrix} \partial_x\, (\coeff \,\partial_x u) = \begin{pmatrix} {\mathel u}_P \\ {\mathel u}_Q\end{pmatrix}\label{pqfull1d}
\end{equation}
where ${\mathel u}_P$ and ${\mathel u}_Q$ denote the low and high-frequency parts of ${\mathel u}(k)$. 
Because we are interested in homogenising $\coeff$, it is written as a matrix, i.e. $\coeff(x) = {\rm diag}(\coeff)$; additionally, the product $\coeff(x)\,u(x)$ is represented as a product
between the diagonal matrix $\coeff$ and vector $u$. We rewrite equation~(\ref{pqfull1d}) as
\begin{eqnarray}
&&\begin{pmatrix} P \\ Q\end{pmatrix}\partial_x\left[\begin{pmatrix} P^* & Q^*\end{pmatrix} \begin{pmatrix} P \\ Q\end{pmatrix} \coeff \begin{pmatrix} P^* & Q^*\end{pmatrix}
\begin{pmatrix} P \\ Q\end{pmatrix}\partial_x\left\{ \begin{pmatrix} P^* & Q^*\end{pmatrix} \begin{pmatrix} P \\ Q\end{pmatrix} u\right\}\right] \nonumber\\
&&\mbox{}= \begin{pmatrix} K_P & 0 \\ 0 & K_Q\end{pmatrix}\begin{pmatrix} P\coeff P^* & P\coeff Q^*\\ Q\coeff P^* & Q\coeff Q^*\end{pmatrix}\begin{pmatrix} K_P & 0 \\ 0 & K_Q\end{pmatrix}\begin{pmatrix} u_P \\ u_Q\end{pmatrix}\nonumber\\
&& \begin{pmatrix} {\mathel u}_P \\ {\mathel u}_Q\end{pmatrix} = \begin{pmatrix} K_P\,P\coeff\,P^*\,K_P & K_P\,P\coeff\,Q^*\,K_Q\\ K_Q\,Q\coeff\,P^*\,K_P & K_Q\,Q\coeff\,Q^*\,K_Q\end{pmatrix}\begin{pmatrix} u_P \\ u_Q\end{pmatrix}.
\end{eqnarray}
 Since we want to model just the coarse scales, we set to zero the fine-scale equation,
 given by the multiplication with the second row, i.e. $\mathel u_Q = 0$, arriving at the equation 
 \begin{eqnarray}
K_Q\, Q\coeff\,Q^*\,K_Q\,u_Q &=& - K_Q\,Q\coeff\,P^*\,K_P\, u_P,\nonumber\\
 u_Q &=&  -(K_Q\, Q\coeff\,Q^*\,K_Q)^{-1}\,(K_Q\,Q\coeff\,P^*\,K_P\, u_P),\label{uq1d}\\
u_Q &=& -K^{-1}_Q( Q\coeff\,Q^*)^{-1}\,Q\coeff\,P^*\,K_P\, u_P.\nonumber 
 \end{eqnarray}
 The coarse-scale part of the equation, given by the multiplication of the first row of the operator with the Fourier projection of $u$,
 \begin{eqnarray}
{\mathel u}_P = K_P[P\coeff\,P^*-P\coeff\,Q^*(Q\coeff\,Q^*)^{-1}\,Q\coeff\,P^*]\,K_P\,Pu,
 \end{eqnarray}
 obtained by substituting the expression for $u^Q$ from equation~(\ref{uq1d}).
Substituting $K_P = P\,\partial_x\,P^*$, and denoting
 \begin{equation}{\bar \coeff} = {\rm diag}(P^*P[\coeff\,-\coeff\,Q^*(Q\coeff\,Q^*)^{-1}\,Q\coeff]P^*P),\label{homog.1d.a} \end{equation}
where the diag operation involves extracting only the terms along the principal diagonal. { Note that we may extract non-diagonal terms from the matrix as well but this corresponds to an integral over space - and thus the original differential equation would become an integro-differential equation. This is explained by the uncertainty principle that governs convolutions and products, described in Section~\ref{uncertainty}.}

The first term, $P^*PaP^*P$ is simply the low-pass representation of $a$ (akin to simply discarding the high wavenumbers) while the second term, the corrector, captures the impact of the high wavenumbers on the low. 
In general the matrix thus obtained will have off-diagonal terms, and one may discard these, especially if it is diagonally dominated. Needless to say,
ignoring the off-diagonal terms will lead to numerical errors when computing solutions. 
Reverting to the original equation, we have,
 \begin{eqnarray}
\mathel u_P &=& P\,\partial_x[ \,{\bar{\coeff}} \,\partial_x {\bar u}],
 \end{eqnarray}
 where ${\bar u} = P^*P u$ is the low-resolution solution, because 
 \begin{eqnarray}
{\bar u}(x) = {\mathcal F}^{-1} \begin{pmatrix} u_P \\ {\bf 0} \end{pmatrix} = \begin{pmatrix} P^* & Q^* \end{pmatrix} \begin{pmatrix} u_P \\ {\bf 0} \end{pmatrix} = P^* P u.\label{pstarp}
 \end{eqnarray}
The high wavenumbers, of which there are $k_q$, have been set to zero in equation~(\ref{uq1d}) which is why we replace $u_Q$ by a vector of zeroes of length $k_q\times1$ in equation~(\ref{pstarp}), and the Fourier coefficient vector is $(u_p\,\, {\bf 0})^T$. 

\begin{figure}[!ht]
\centering
\epsfig{file=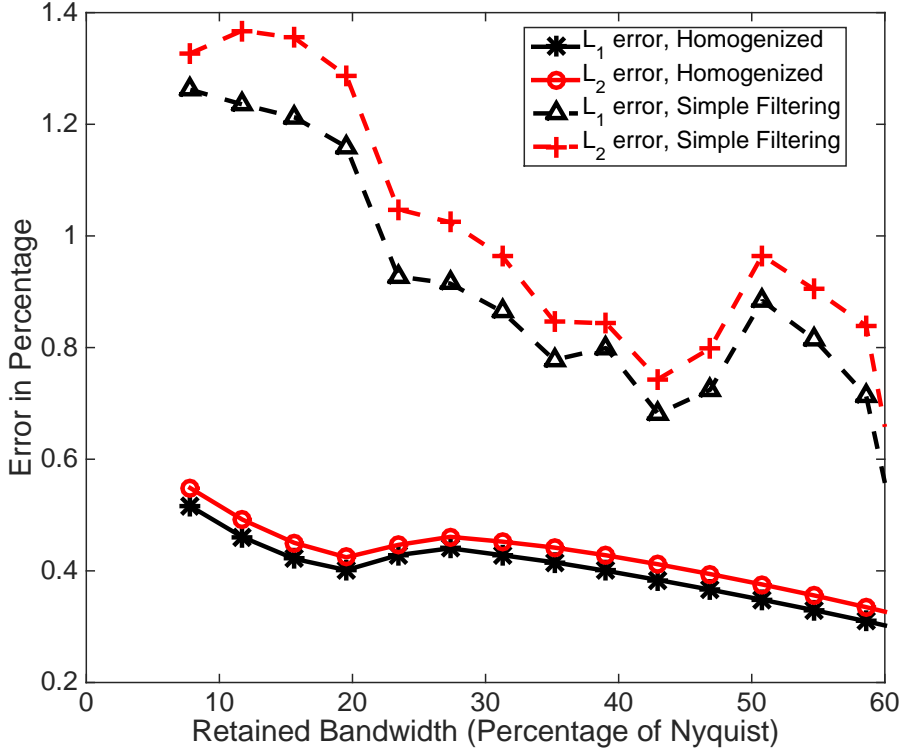,width=\linewidth}
\caption{Homogenization consistently outperforms raw filtering in 1D: errors in decimating heterogeneous coefficients of the diffusion equation, $\mathel u = \partial_x (\coeff\partial_x u) = 0$, plotted against the amount of bandwidth retained. For instance, 0 retention corresponds to replacing the heterogeneous coefficient by one number and 100 means retaining the original coefficient. In order to satisfy Orszag's two-thirds rule \citep{orszag71}, we zero out the upper third portion of the spectrum, which is why the error in the raw filtered coefficient rapidly drops to zero around 66\% of the Nyquist. The $x$-axis range therefore spans 0 to 60\%.
Setting boundary conditions on the domain $x\in [0,1], u(0) = 0, u(1) = 1$, the exact solution can be obtained: $u(x) = \left(\int_0^1 \coeff^{-1}(x') dx'\right)^{-1}\,\int_0^x \coeff^{-1}(x') dx' $. 
The homogenised and raw filtered coefficients are computed, and the difference errors measured in the $L_1$ and $L_2$ norms, normalised by the corresponding norms of the exact solution, are plotted.  It is however important to note that these curves are sensitive to the model of the coefficient $\coeff$ and the curves should be thought of more as indicative of the trend rather than an absolute statement on the errors. Although the problem is not periodic, the effective coefficients are still estimated in the Fourier basis.
\label{1dhomogenize}}
\end{figure}

A first test is to compare the outcome of equation~(\ref{homog.1d.a}) to the classical homogenisation result. For a periodic $\coeff$ in the diffusion equation $\mathel u = \partial_x (\coeff\partial_x u) = 0$, the homogenised coefficient is $[\int_0^1\, dx'\, \coeff^{-1}(x')]^{-1}$, with $x\in[0,1]$, a result that we obtained successfully \citep[the derivation is fully laid out for instance in][]{varadhan82,capdeville1d}. 
Next considered the 1D diffusion equation with a heterogeneous coefficient. The exact solution to the problem in $x\in [0,1], u(0) = 0, u(1) = 1$ is given by $u(x) = \left(\int_0^1 \coeff^{-1}(x') dx'\right)^{-1}\,\int_0^x \coeff^{-1}(x') dx' $. We compared the exact solution to that obtained using the homogenized and raw-filtered coefficients. The error is plotted as a function of retained bandwidth of the coefficient in Figure~\ref{1dhomogenize} and shows that homogenisation is distinctly superior to there performance of raw filtering.

\subsection{RG notation}
The development outlined in section~(\ref{num1d}) can be recast in the notation of RG. The Fourier domain is naturally suited for RG; we use subscripts `$<$' and `$>$' to describe low and high wavenumbers and the goal is to remove all the $>$ (high) wavenumbers, projecting them onto the $<$ (low) wavenumbers. We recall that `high' and `low' are defined externally, suitable for the problem at hand. { The continuous Fourier transform of $\mathel u$ is 
\begin{equation}
{\mathel u}(k) = - \int_{-\infty}^\infty dk' \,k\,k'\,\coeff(k-k')\,u(k').\label{kkeq}
\end{equation}
Denoting wavenumbers $|k| \le k_c$ by $k_<$ and wavenumbers $|k| > k_c$ by $k_>$, where $k_c$ is the desired cutoff, we define the following integrals
\begin{eqnarray}
\int_< dk  = \int_{-k_c}^{k_c} dk, \,\,\,\,\,\,\, \int_> dk  = \int_{-\infty}^{-k_c} dk + \int_{k_c}^{\infty} dk .\label{eqintdef}
\end{eqnarray}
Using this notation, we can split equation~(\ref{kkeq}) into two parts, for small and high $k$ respectively,
\begin{eqnarray}
{\mathel u}(k)_< = - \int_{<} dk' \,k_<\,k'\,\coeff(k-k')\, u_<(k') - \int_{>} dk \,k_<\,k'\,\coeff(k-k')\, u_>(k'),\label{low.rg}\\
{\mathel u}(k)_> = - \int_{>} dk' \,k_>\,k'\,\coeff(k-k')\, u_>(k')- \int_{<} dk' \,k_>\,k'\,\coeff(k-k')\, u_<(k').\label{high.rg}
\end{eqnarray}
Equation~(\ref{low.rg}) for instance states that $k$ belongs to the set of low wavenumbers, and therefore if $k'$ is also a low wavenumber,
$ u(k') = u_<$ and if $k'$ were large, $ u(k') =  u_>$. As described in equation~(\ref{eqintdef}), integrals with the subscript $<$ represent integration over low wavenumbers whereas
those with the subscript $>$ denote high-wavenumber integration.
%
The idea thus involves integrating over shells of the high wavenumbers and expressing them in terms of the low wavenumbers,
\begin{equation}
{\mathel u}(k)_> = 0 \rightarrow  \int_{>} dk' \,k_>\,k'\,\coeff(k-k')\, u_>(k') = - \int_{<} dk' \,k_>\,k'\,\coeff(k-k')\, u_<(k').\label{projections}
\end{equation}
Equation~(\ref{projections}) details the connection between $u_>$ and $u_<$ through a Fredholm integral over the range $k'\in H$. Positing the existence an inverse Fredholm operator $a^{-1}(k,k')$,
we write,
\begin{equation}
{u}(k)_> = -\int_> dk'\,\int_< dk{''}\,\coeff^{-1}(k_>,k')\, \coeff(k',k^{''})\,{u}_<(k^{''}).
\end{equation}
In general, the real-space representation of the high-frequency portion of the convolution operator is $\coeff(k_>,k'_>) =  \coeff_{>>}(x,x')$, where $x,x'$ are spatial variables corresponding to Fourier wavenumbers $k,k'$.
Whereas the Fourier transform of the coefficient $\coeff$ contains only one wavenumber, the convolutional integrals in, e.g. equation~(\ref{low.rg}), introduce two wavenumbers $k,k'$ and thus the originally one-variable function $\coeff$ may become a function of two spatial coordinates when split into high- and low-wavenumber combinations $<<, <>, ><, >>$.
\subsection{Uncertainty principle applied to convolutions}\label{uncertainty}
A deeper connection may be found when considering the coefficient using the theory of distributions. The product of two functions can be written as integral, i.e.
\begin{equation}
a(x)\,u(x) = \int_{-\infty}^\infty dx'\,a(x')\delta(x-x')\,u(x'),
\end{equation}
and so we consider the redefinition of the coefficient $a$ as an operator, 
\begin{equation}
\coeff(x,x') =  \coeff(x')\, \delta(x-x') =\int_{-\infty}^\infty dk'\,\int_{-\infty}^\infty dk\,\hat{\coeff}(k')\, e^{ik(x-x') + ik'x'},\label{delta.op}
\end{equation}
and substituting the following transformation $\eta = k - k'$ and replacing $k'$ by $k-\eta$, we obtain
\begin{eqnarray}
&&\int_{-\infty}^\infty d\eta\,\int_{-\infty}^\infty dk\,\hat{\coeff}(k-\eta)\, e^{-i\eta x' + ik x} = \\
&&\left[\int_{<} d\eta\,\int_< dk + \int_{<}d\eta\, \int_> dk + \int_{>}d\eta\, \int_< dk + \int_{>}d\eta\,\int_> dk \right] \,[\hat{\coeff}(k-\eta)\, e^{-i\eta x' + ik x}].\nonumber
\end{eqnarray}
It can be seen that only if all four integrals are retained, the resultant
quantity reduces to the form in equation~(\ref{delta.op}). If however, one of the terms, such as the $> >$ integral, is projected on to the others, the inverse transform produces
a more general function $\bar{a}(x,x')$.} 
It is possible to connect the degree of broadening in space (i.e. the deviation of the operator $a(x,x')$ from $a(x)\delta(x-x')$) to the degree of squeezing in Fourier domain, i.e. the number of Fourier coefficients that are discarded. Indeed, \cite{stegel00} has derived an
uncertainty principle for convolution operators, which address the tradeoff between spectral and real-space representations. In addition, the Fredholm integral equation~(\ref{projections}) typically appears in the classical {\it spectral concentration problem} \citep{slepian} that concerns the design of functions that are optimally spatio-spectrally localized. The analogy with the spatio-spectral concentration of convolution operators thus follows.

\begin{figure}[!ht]
\centering
\epsfig{file=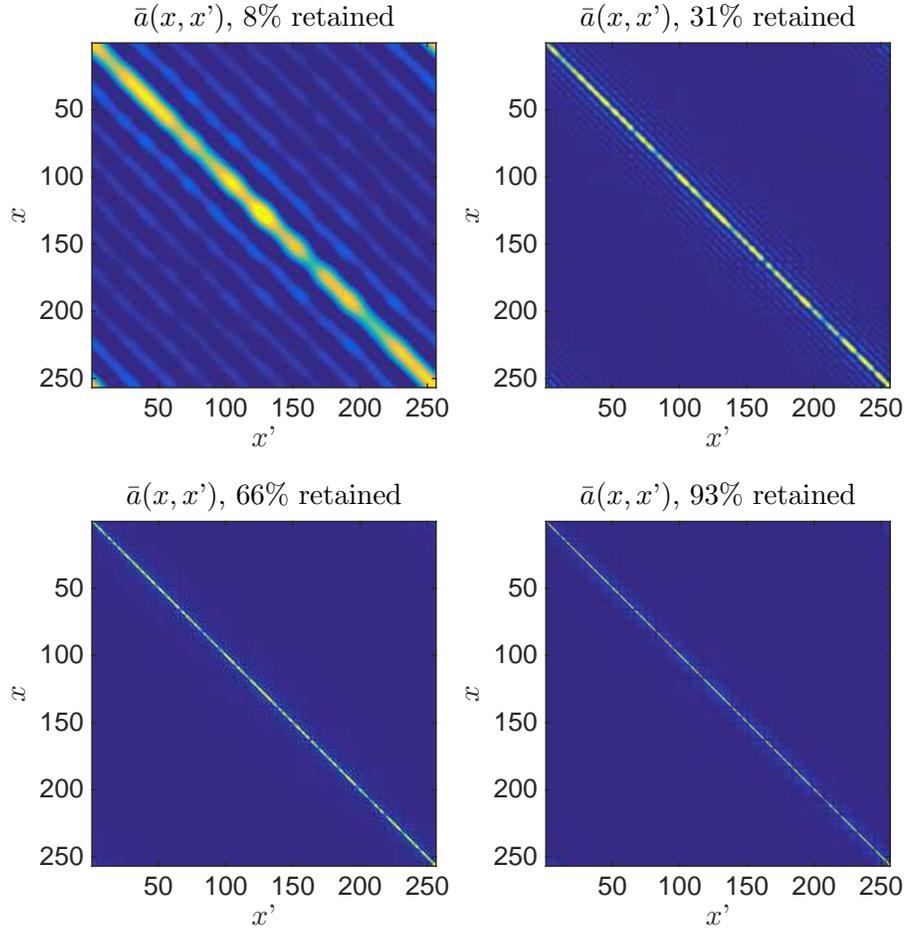,width=\linewidth}
\caption{Uncertainty principle for convolutions: spatial broadening as a consequence of Fourier-domain squeezing. As discussed in section~\ref{num1d}, homogenisation contributes to the emergence of finite off-diagonal terms in $\bar{a}(x,x')$. As smaller fractions of the bandwidth are retained, the broadening in $(x,x')$ space increases. The colour scale is the same on all plots. A Heisenberg uncertainty principle operates; squeezing in Fourier domain results in broadening in physical space \citep{stegel00}. However, unlike the classical theorem that applies to transform pairs, this is a convolutional uncertainty, in that squeezing a convolution in the $(k,k')$ space causes the $(x,x')$ representation to broaden.
\label{xxp}}
\end{figure}

In linear equations, especially ones that
only involve space in the diffusion equation, ignoring the broadening in real space can be a source of numerical error.
Substituting back into the inverse Fourier transformed form of equation~(\ref{low.rg}),
\begin{eqnarray}
\mathel u_< &=& \partial_x(\bar{\coeff}\,\partial_x{u}_<),\\
\bar{\coeff}(x,x') &=& \int_< dk \int_< dk'\, e^{ikx + ik' x'} \,\left[\coeff(k,k')\right.\\
 &-&\left.\int_> dk^{''}\int_> dk^{'''}\,  \, \coeff(k,k^{''})\,\coeff^{-1}(k^{''},k^{'''})\, \coeff(k^{'''},k')\right],\nonumber
\end{eqnarray}
which is essentially what equation~(\ref{homog.1d.a}) states but in continuous notation. { As in equation~(\ref{homog.1d.a}), the renormalised coefficient is given by $\bar{\coeff} = \bar{\coeff}(x,x')$ and the operator acts on a function $u(x)$ thus, 
\begin{equation}
\bar{\coeff}u = \int dx'\, {\bar a}(x,x')\, u(x').
\end{equation}
This introduces a greater degree of computational complexity, in the sense that what was originally a differential equation is now an integro-differential equation. As a consequence, we force $\bar{\coeff} \rightarrow \bar{\coeff}(x,x')\,\delta(x-x')$, potentially leading to errors in the coarse-scale representation of $\mathel u$.} Figure~\ref{xxp} illustrates the manifestation of the uncertainty principle. As the retained bandwidth of the convolution operator reduces, so expands its representation in real space, evidenced by the off-diagonal components. { Note that the approximation $\bar{\coeff} \rightarrow \bar{\coeff}(x,x')\,\delta(x-x')$ translates to dropping the off-diagonal terms (see Figure~\ref{xxp}).}
 

%

%
%

\section{2D case}
We consider the diffusion equation with a scalar coefficient in 2D:
\begin{equation}
\bnabla\cdot(\coeff\bnabla u) = 0,\label{system1}
\end{equation}
where $\bnabla$ is the covariant spatial gradient, $\coeff = \coeff(\bx) > 0$, $\bx$ is the 2D Cartesian coordinate $\bx = (x,y)$ and $u=u(\bx)$ is the solution. 
Fourier transforms in the $x$ and $y$ directions are denoted by $\ft_x = \begin{pmatrix} P_x \\ Q_x\end{pmatrix}$ and $\ft_y = \begin{pmatrix} P_y \\ Q_y\end{pmatrix}$.  The 2D transform commutes in that $\ft_x \ft_y = \ft_y \ft_x$ and is written thus
\begin{eqnarray}
\begin{pmatrix} P_x P_y \\ Q_x P_y \\ P_x Q_y\\ Q_x Q_y\end{pmatrix},
\end{eqnarray}
where the first entry is the low-pass filter and the three other entries represent high-pass filters. 
The following identities hold for pairwise products of the projection matrices
\begin{eqnarray}
&& P_x P_y = P_y P_x,\,\,\, Q_x P_y = P_y Q_x,\,\,\,P_x Q_y = Q_y P_x, \,\,\, Q_x Q_y = Q_y Q_x,\label{2dids}\\
&&P_xP_x^* = \eye,\,\,\, Q_xQ_x^* = \eye,\,\,\, P_xQ_x^* = {\bf 0},\,\,\,
Q_xP_x^* = {\bf 0},\,\,\, P_x^*P_x + Q_x^*Q_x = \eye,\nonumber\\
&&P_yP_y^* = \eye,\,\,\, Q_yQ_y^* = \eye,\,\,\, P_yQ_y^* = {\bf 0},\,\,\,
Q_yP_y^* = {\bf 0},\,\,\, P_y^*P_y + Q_y^*Q_y = \eye.\nonumber
\end{eqnarray}
{ Because the $y$-projection matrices are unaffected by the operator $\partial_x$ and vice versa, we have similar relationships as in equation~(\ref{2dids}). We do not explicitly write them here.}

The projection of the partial spatial derivative $\partial_x$ on to 2D Fourier space is 
\begin{eqnarray}
\begin{pmatrix} P_x P_y \\ Q_x P_y \\ P_x Q_y\\ Q_x Q_y\end{pmatrix} \partial_x \begin{pmatrix} P_x^* P_y^* & Q_x^* P_y^* & P_x^* Q_y^* & Q_x^* Q_y^*\end{pmatrix}
 = \begin{pmatrix} K^P_x & 0 & 0 & 0 \\ 0 & K^Q_x & 0 & 0\\ 0 & 0  &K^P_x & 0 \\ 0 & 0 & 0 &K^Q_x\end{pmatrix},\nonumber
\end{eqnarray}
through the use of identities~(\ref{2dids}), which we write in anticipation of future notational needs as 
\begin{equation}
\begin{pmatrix} K^P_x & 0 & 0 & 0 \\ 0 & K^Q_x & 0 & 0\\ 0 & 0  &K^P_x & 0 \\ 0 & 0 & 0 &K^Q_x\end{pmatrix} =  \begin{pmatrix} K^P_x & 0\\ 0 & \mathk_x \end{pmatrix}.\label{defs2dx}
\end{equation}

The projection of the partial spatial derivative $\partial_y$ is 
\begin{eqnarray}
\begin{pmatrix} P_x P_y \\ Q_x P_y \\ P_x Q_y\\ Q_x Q_y\end{pmatrix} \partial_y \begin{pmatrix} P_x^* P_y^* & Q_x^* P_y^* & P_x^* Q_y^* & Q_x^* Q_y^*\end{pmatrix}
 = \begin{pmatrix} K^P_y & 0 & 0 & 0 \\  0& K^P_y & 0 & 0\\ 0 & 0  &K^Q_y & 0 \\ 0 & 0 & 0 &K^Q_y\end{pmatrix},
\end{eqnarray}
rewritten as
\begin{equation}
\begin{pmatrix} K^P_y & 0 & 0 & 0 \\  0& K^P_y & 0 & 0\\ 0 & 0  &K^Q_y & 0 \\ 0 & 0 & 0 &K^Q_y\end{pmatrix}
=  \begin{pmatrix} K^P_y & 0\\ 0 & \mathk_y \end{pmatrix}.\label{defs2dy}
\end{equation}

The projection of $\coeff$ on to this basis is
\begin{equation}
\begin{pmatrix} P_x P_y \\ Q_x P_y \\ P_x Q_y\\ Q_x Q_y\end{pmatrix}\coeff \begin{pmatrix} P_x^* P_y^* & Q_x^* P_y^* & P_x^* Q_y^* & Q_x^* Q_y^*\end{pmatrix} = \begin{pmatrix} D & C \\ B & A\end{pmatrix},
\end{equation}
where
\begin{eqnarray}
&&D =  P_xP_y\,\coeff\,P^*_xP^*_y,\nonumber\\
&& C = \begin{pmatrix} P_xP_y\,\coeff\,Q^*_xP^*_y & P_xP_y\,\coeff\,P^*_xQ_y^* & P_xP_y\,\coeff\,Q^*_xQ_y^* \end{pmatrix},\nonumber\\ 
&& B = \begin{pmatrix} Q_xP_y\,\coeff\,P^*_xP^*_y \\ P_xQ_y\,\coeff\,P^*_xP^*_y \\ Q_xQ_y\,\coeff\,P^*_xP^*_y \end{pmatrix},\nonumber\\ 
&& A = \begin{pmatrix}Q_xP_y\,\coeff\,Q^*_xP^*_y & Q_xP_y\,\coeff\,P^*_xQ^*_y & Q_xP_y\,\coeff\,Q^*_xQ^*_y
\\ P_xQ_y\,\coeff\,Q^*_xP^*_y  &P_xQ_y\,\coeff\,P^*_xQ^*_y & P_xQ_y\,\coeff\,Q^*_xQ^*_y \\ 
Q_xQ_y\,\coeff\,Q^*_xP^*_y & Q_xQ_y\,\coeff\,P^*_xQ^*_y &Q_xQ_y\,\coeff\,Q^*_xQ^*_y\end{pmatrix}.
\end{eqnarray}
The projection of the solution on to Fourier space is written as
\begin{equation}
\begin{pmatrix} P_x P_y \\ Q_x P_y \\ P_x Q_y\\ Q_x Q_y\end{pmatrix}u = \begin{pmatrix} u_P \\ u_Q \end{pmatrix},
\end{equation}
where $u_P = P_x P_y\,u$ is the low-pass filtered component of the solution
and $u_Q$ contains the high-frequency projections in the $x$ and $y$ directions. 
The full projection is 
\begin{equation}
\begin{pmatrix} K^P_x & 0\\ 0 & \mathk_x \end{pmatrix}\begin{pmatrix} D & C \\ B & A\end{pmatrix}\begin{pmatrix} K^P_x & 0\\ 0 & \mathk_x \end{pmatrix}\begin{pmatrix} u_P \\ u_Q\end{pmatrix}
 + \begin{pmatrix} K^P_y & 0\\ 0 & \mathk_y \end{pmatrix}\begin{pmatrix} D & C \\ B & A\end{pmatrix}\begin{pmatrix} K^P_y & 0\\ 0 & \mathk_y \end{pmatrix} \begin{pmatrix} u_P \\ u_Q\end{pmatrix}\nonumber,
\end{equation}
or
\begin{eqnarray}
=\begin{pmatrix} K^P_xDK^P_x & K^P_xC\mathk_x \\ \mathk_x B K^P_x & \mathk_x A\mathk_x\end{pmatrix}\begin{pmatrix}u_P \\ u_Q\end{pmatrix} +
\begin{pmatrix} K^P_yDK^P_y & K^P_yC\mathk_y \\ \mathk_y B K^P_y & \mathk_y A\mathk_y\end{pmatrix}\begin{pmatrix}u_P \\ u_Q\end{pmatrix},\nonumber\\
=\begin{pmatrix} K^P_xDK^P_x + K^P_yDK^P_y & K^P_xC\mathk_x + K^P_yC\mathk_y \\ \mathk_x B K^P_x + \mathk_y B K^P_y & \mathk_x A\mathk_x + \mathk_y A\mathk_y\end{pmatrix}\begin{pmatrix}u_P \\ u_Q\end{pmatrix}.
\end{eqnarray}
As in section~\ref{num1d}, setting the high-frequency component to zero, we obtain
\begin{equation}
u_Q = - (\mathk_x A\mathk_x + \mathk_y A\mathk_y)^{-1} (\mathk_x B K^P_x + \mathk_y B K^P_y)u_P,
\end{equation}
and substituting this into the low-frequency component, we obtain
\begin{eqnarray}
[(K^P_xDK^P_x + K^P_yDK^P_y) - \nonumber\\
(K^P_xC\mathk_x + K^P_yC\mathk_y)(\mathk_x A\mathk_x + \mathk_y A\mathk_y)^{-1}(\mathk_x B K^P_x + \mathk_y B K^P_y)] u_P = 0,\nonumber
\end{eqnarray}
and recalling the definitions of $K^P_x$ and $K^P_y$ from equations~(\ref{defs2dx}) and~(\ref{defs2dy}), we arrive at the effective low-pass equation
\begin{equation}
\bnabla\cdot(\bar{\bcoeff}\cdot\bnabla{\bar u}) = 0,
\end{equation}
where $\bar{\bcoeff}$ is a second-order tensor
\begin{equation}
\bar{\bcoeff} = \begin{pmatrix}\bar\coeff_{xx} &\bar\coeff_{xy} \\ \bar\coeff_{yx} & \bar\coeff_{yy}\end{pmatrix}.\label{tensor.coeff}
\end{equation}
The different components of the coefficients are given by,
\begin{eqnarray}
\bar\coeff_{xx} = {\rm diag}\{P_x^*P_y^*[D -  C\,\mathk_x(\mathk_x A\mathk_x + \mathk_y A\mathk_y)^{-1}\mathk_x\, B]P_xP_y\},\\
\bar\coeff_{xy} = -{\rm diag}\{P_x^*P_y^*\,[C\,\mathk_x(\mathk_x A\mathk_x + \mathk_y A\mathk_y)^{-1}\mathk_y\, B]\,P_xP_y\},\\
\bar\coeff_{yx} = -{\rm diag}\{P_x^*P_y^*\,[C\,\mathk_y(\mathk_x A\mathk_x + \mathk_y A\mathk_y)^{-1}\mathk_x\, B]\,P_xP_y\},\\
\bar\coeff_{yy} = {\rm diag}\{P_x^*P_y^*[D -  C\,\mathk_y(\mathk_x A\mathk_x + \mathk_y A\mathk_y)^{-1}\mathk_y\, B]P_xP_y\}.
\end{eqnarray}
With no conceptual difficulty (but a greater degree of book keeping), this formulation may be extended to equations with tensor coefficients, outlined for instance in Appendix~\ref{tensor.eq}. 

\begin{figure}[!ht]
\centering
\epsfig{file=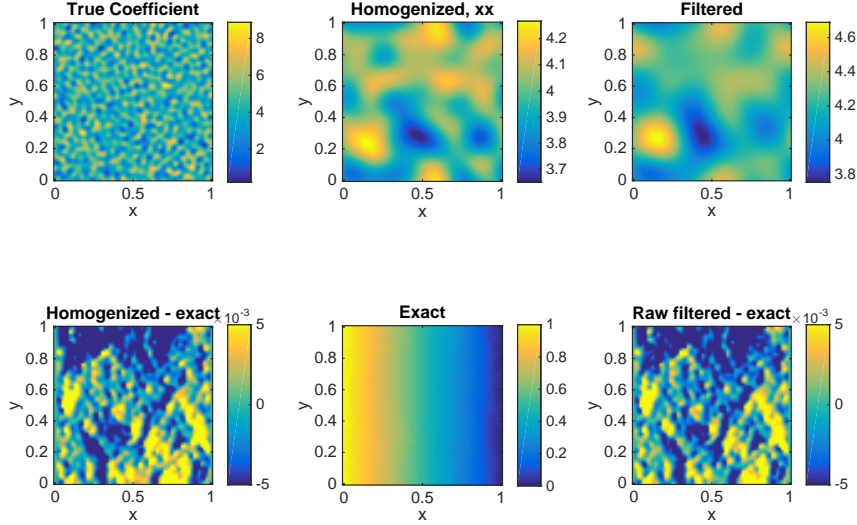,width=\linewidth}
\caption{Finite-wavenumber homogenisation applied to the 2D diffusion equation with a scalar coefficient (upper-left panel). 
The homogenised ($\bar \coeff_{xx}$) and (raw) filtered coefficients, shown in the upper-middle and upper-right panels.
shows the exact solution, the lower-left panel the difference between the homogenised and exact and the lower-right, the difference between the raw-filtered and exact solutions. Only 2 coefficients are retained from the full solution. There is little perceivable difference between the homogenised and filtered solutions.{ We use the freely available Portable Extensible Toolkit for Scientific Computation (PETSc) to compute the solution on a $64\times64$-sized grid applying second-order centered finite differences. Zero-Neumann boundary conditions are applied to the upper and lower boundaries, $\partial_y u(x,y=0) = 0 = \partial_y u(x,y=1)$, unit-Dirichlet on the left, $u(x=0,y) = 1$ and zero-Dirichlet on the right, $u(x=1,y) = 0$ boundaries.} Note that as in Figure~\ref{1dhomogenize}, we are careful to honour Orszag's two-thirds rule \citep{orszag71}.
\label{2dhomogenize}}
\end{figure}

The first and simplest test of this method is to homogenise a periodic coefficient that varies only in one direction (e.g. \cite{engquist02}), i.e. where $\coeff = \coeff(x)$, and whose periodicity length scale is small in comparison to the overall size of the domain. Classical homogenisation predicts a splitting of the scalar coefficient into tensorial components, $\bar{\coeff}_{xx} = (\int_0^1 dx' \coeff^{-1}(x'))^{-1}$ and $\bar{\coeff}_{yy} = (\int_0^1 dx' \coeff(x'))$. We solve the 2D scalar-coefficient diffusion equation, shown in Figure~\ref{2dhomogenize}, where we build a 2D model of the coefficient $\coeff$ and decimate it using methods of homogenisation and raw filtering. { The coefficients obey $\coeff_{xx} = \coeff_{yy}, \coeff_{xy} = 0 = \coeff_{yx}$. We first populate all the spatial points of the coefficient using a uniform random distribution and subsequently filter out the upper third of the Fourier transform (i.e. by setting it to zero) to conform with Orszag's two-thirds rule \citep{orszag71}. A water level is added to the coefficient to ensure that its minimum is finite and positive (see Figure~\ref{2dhomogenize}).}

\begin{figure}[!ht]
\centering
\epsfig{file=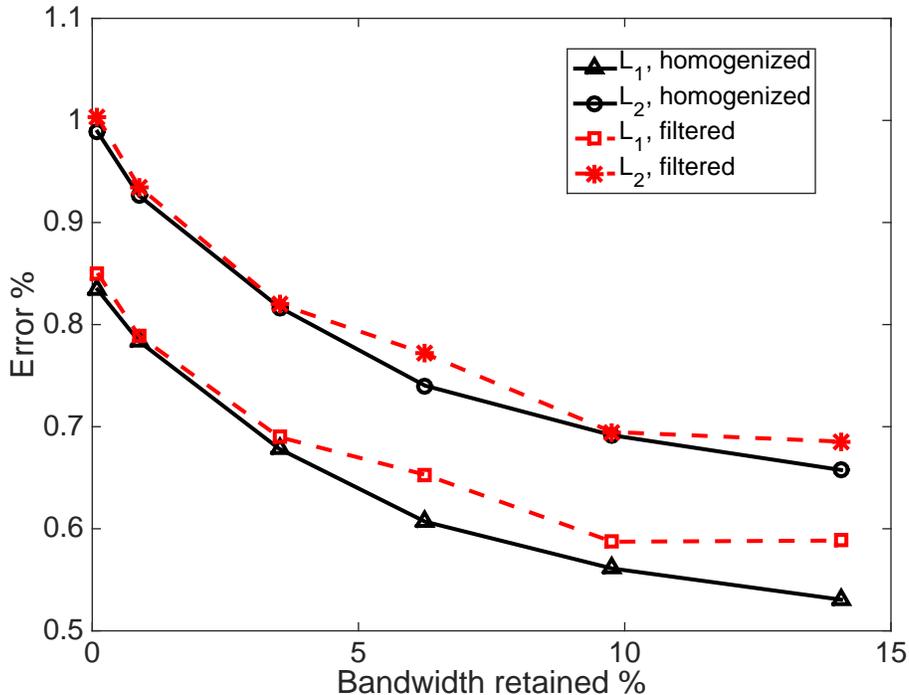,width=\linewidth}
\caption{Homogenization appears to perform better or approximately as well as raw filtering in 2D: errors as measured in the $L_1, L_2$ norms between the exact and the homogenised and filtered solutions as a function of the retained bandwidth of the coefficient for the case shown in Figure~\ref{2dhomogenize}. We solve the diffusion problem as described in Figure~\ref{2dhomogenize}. Note that homogenisation produces a tensor coefficient whereas raw filtering produces a low-pass model of the original coefficient, i.e. $\coeff_{xx} = \coeff_{yy} = {\rm filt}(\coeff)$. Homogenisation does better at intermediate bandwidths, i.e. between zero wavenumber and at a third of the Nyquist (Orszag's two-thirds rule). The coefficient has no power beyond  a third of the Nyquist, and therefore filtering is seen to outperform (or approximately as well as) homogenisation.
\label{xxp2}}
\end{figure}

Homogenization produces a tensor coefficient at coarse scales, as in equation~(\ref{tensor.coeff}), whereas raw filtering just removes the high spatial frequencies of the scalar coefficient. Comparing solutions obtained using these two decimated-coefficient models in Figure~\ref{xxp2}, we find that the homogenized solution performs slightly better than the naively filtered case. Homogenisation does not perform significantly better than raw filtering, as was observed in the 1-D case, and the reason is not apparent. 

\begin{figure}[!ht]
\centering
\epsfig{file=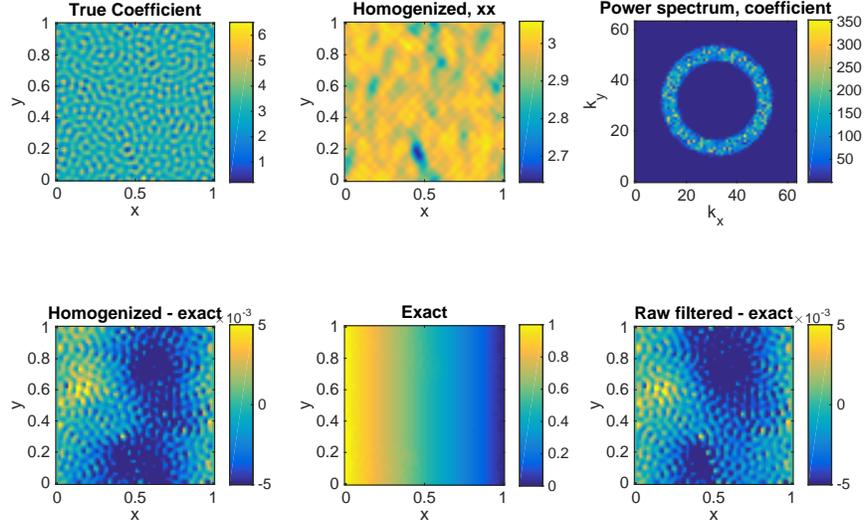,width=\linewidth}
\caption{Finite-wavenumber homogenisation applied to the 2D diffusion equation with a scalar coefficient (upper-left panel) that has a sparse Fourier representation (upper-right panel). 
Nominally, the coefficient would have had power up to a third of the Nyquist (honouring Orszag's two-thirds rule) but to enforce sparsity we limit the wavenumber range over which the coefficient has power. Note that in order to improve visibility, we have subtracted the mean from the coefficient (zero wavenumber power) in the upper-right panel. The homogenised ($\bar \coeff_{xx}$) coefficient is shown in the upper-middle panel.
The lower-middle panel shows the exact solution, the lower-left panel the difference between the homogenised and exact and the lower-right, the difference between the raw-filtered and exact solutions.{ We use the publicly downloadable Portable Extensible Toolkit for Scientific Computation (PETSc) to compute the solution on a $64\times64$-sized grid applying second-order centered finite differences. Zero-Neumann boundary conditions are applied to the upper and lower boundaries, $\partial_y u(x,y=0) = 0 = \partial_y u(x,y=1)$, unit-Dirichlet on the left, $u(x=0,y) = 1$ and zero-Dirichlet on the right, $u(x=1,y) = 0$ boundaries.} Note that as in Figure~\ref{1dhomogenize}, we are careful to honour Orszag's two-thirds rule \citep{orszag71}. 
\label{sparsehom}}
\end{figure}

\begin{figure}[!ht]
\centering
\epsfig{file=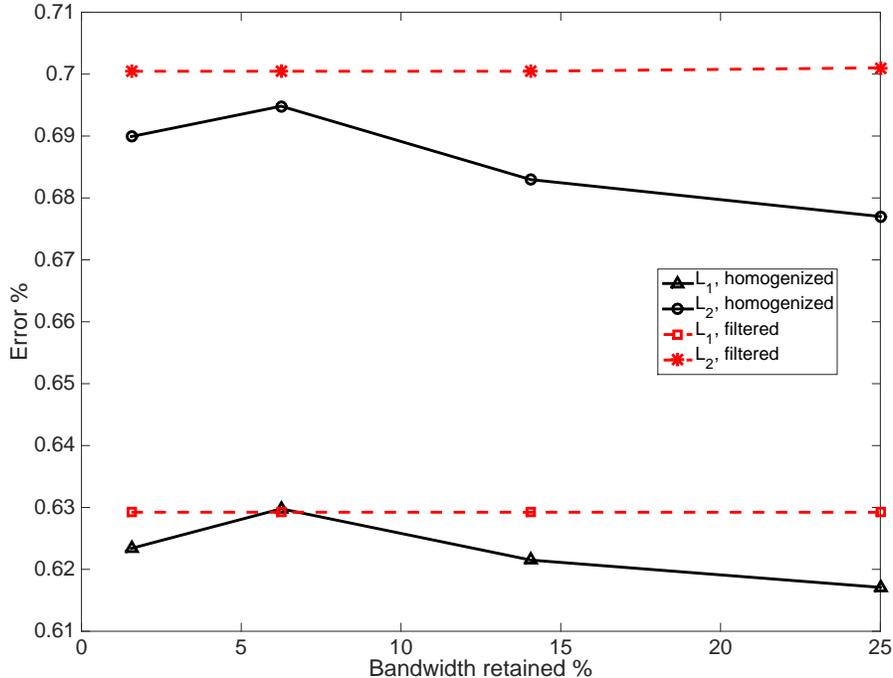,width=\linewidth}
\caption{Homogenisation performs better than raw filtering when the coefficient is sparse in the basis (Figure~\ref{sparsehom}): errors as measured in the $L_1, L_2$ norms between the exact and the homogenised and filtered solutions as a function of the retained bandwidth of the coefficient. We solve the diffusion problem as described in Figure~\ref{2dhomogenize}. Note that homogenisation produces a tensor coefficient whereas raw filtering produces a low-pass model of the original coefficient, i.e. $\coeff_{xx} = \coeff_{yy} = {\rm filt}(\coeff)$. Because there is no power between zero wavenumber and the narrow annulus where all the power resides, the errors associated raw filtering are relatively flat over a range. When the filter begins to include the annulus, the raw filtering error starts to drop sharply.
\label{error.sparse}}
\end{figure}

 Finally we study the effect of basis sparsity of the coefficient has on the quality of the homogenised solution. We generate coefficients that are sparse in the Fourier basis (Figure~\ref{sparsehom}) and compute solutions. The errors between the exact and homogenised / filtered solutions as a function of retained spectral bandwidth are plotted in Figure~\ref{error.sparse}. Note that because the power is narrowly focused, raw filtering errors are roughly constant to the point when the filter begins to include the annulus (see upper-right panel in Figure~\ref{sparsehom}), at which point the error drops sharply. Homogenisation appears to perform well at intermediate ranges of retained bandwidth, away from zero wavenumber and a third of the Nyquist (keeping in mind Orszag's two-thirds rule).

\section{Conclusions}
Numerical homogenization is a powerful methodology to build sub-grid-scale models of coefficients of differential equations
in a variety of scenarios.
The technique of numerical homogenisation bears a striking resemblance to mode-elimination Renormalization Group (RG) 
as developed by Wilson \citep{wilson1975}. Indeed, RG has been applied to a wide range of problems (mostly of the non-linear kind), 
including in fluid mechanics, specifically to model turbulence \citep{yakhot86}. 
In this article, we apply RG to linear problems, viz. the diffusion and wave equations with deterministic coefficients.
We describe how to build sub-grid-scale models (coarse descriptors) of heterogeneous coefficients of differential equations and
characterise the accuracy of the approximation as a function of the degree of coarsening. This raises the more fundamental question
of why RG works well in some scenarios and not others. We trace the problem to the uncertainty principle, which can be derived for
convolutions also \citep{stegel00}. A product of functions in real space is encoded as a convolution in Fourier domain; the extended uncertainty
principle then states that localising convolutions in Fourier domain results in a broadening in real space. In other words, squeezing in $(k,k')$-space
will cause the function to become delocalized in $(x,x')$ space. This observation may be extended to differential equations that evolve temporally in that
squeezing in the spatial Fourier domain will not only broaden the real-space representation but additionally, might introduce temporal convolutions.
Thus RG will likely not work in situations where the convolution is strongly compressed in Fourier domain but where the spatial representation does not
account for delocalisation, i.e. where the coefficients continue are forced to assume the form $\bar\coeff(x,x') = \delta(x-x')\,\coeff(x)$. 
Constructing the right basis on which to project the coefficients can therefore be critical to generating accurate coarsened models. 
Linear problems provide an excellent opportunity to gain insights of this sort into RG and can additionally shed light on why sub-grid-scale models in turbulence
perform well on some occasions but not on others. 

Appendix~\ref{tensor.eq} analyzes the time-dependent wave equation and appears to suggest that the uncertainty principle operates in an unexpected fashion: concentrating the convolution operator in spatial domain results in broadening in both space and {\it time} domains. 
Solutions to the linear wave equation may be projected onto the eigen-basis of the operator, and discarding high spatial wavenumbers implies that the corresponding eigenfunctions are also removed. Because normal modes of the wave operator are at specific wavenumbers and temporal frequencies, selectively removing some of them has consequences in both spatial and temporal domains. This will likely result in broadening in both space and time; the uncertainty principle will thus likely involve the spatio-temporal density of normal modes of the wave equation. 


We have studied relatively simple problems in this article and much work needs to be done
to fortify these results. One of the more serious aspects that remains to be addressed when solving large 3D homogenization 
problems is that of storing and computing the inverses of these large matrices. Sparsity is therefore a critical feature.
The Fourier basis has poor sparsity properties and consequently, wavelet or other bases may prove superior. 
It is also important to characterize the efficiency and accuracy of the method over a wider range of problems
and more realistic scenarios.
\section*{Acknowledgements}
The author acknowledges support from consultancy project PT54324 with Shell India, Ramanujan fellowship SB/S2/RJN-73/2013, the Max-Planck Partner Group Program and the NYU Abu Dhabi Center for Space Science. The work has benefited substantially from numerous conversations, and SMH gratefully acknowledges discussions with Alain Plattner, Srinivasa Varadhan, Rishi Sharma and Sandip Trivedi.

\bibliography{ms.bbl}

\appendix \section{Tensor equation}\label{tensor.eq}
Similar extensions to higher-order wavespeed tensors in 3-D are possible. 
We only outline the method, leaving its execution and testing to a future occasion.
Consider the general seismic wave equation
\begin{eqnarray}
\rho\partial_t^2\bu - \bnabla\cdot({\bf T}:\bnabla\bu) = {\bf f}(\bx, t),
\end{eqnarray}
where ${\bf T} = \{\tau_{ijkl}\}$ is a fourth-order tensor \citep{DT98}. We transform the equation to temporal Fourier domain
\begin{eqnarray}
- \rho\omega^2 \bu - \bnabla\cdot({\bf T}:\bnabla\bu) = {\bf f}(\bx,\omega),
\end{eqnarray}
where for convenience we have not explicitly indicated that $\bu$ is the temporal Fourier trasform.
The projection of the spatial operator on to the Fourier basis is written as
\begin{equation}
\sum_{i,k} \begin{pmatrix}K^P_i & 0\\ 0&\mathk_i\end{pmatrix} \begin{pmatrix}D_{ijk\ell} & C_{ijk\ell}\\ B_{ijk\ell}&A_{ijk\ell}\end{pmatrix}  \begin{pmatrix}K^P_k & 0\\ 0&\mathk_k\end{pmatrix} \begin{pmatrix}u^P_{\ell} \\ u^Q_{\ell}\end{pmatrix},\label{prod}
\end{equation}
where
\begin{equation}
{\mathcal F}_x {\mathcal F}_y {\mathcal F}_z \,\,\tau_{ijk\ell}\,\, {\mathcal F}^{-1}_x {\mathcal F}^{-1}_y {\mathcal F}^{-1}_z
= \begin{pmatrix} D_{ijk\ell} & C_{ijk\ell}\\ B_{ijk\ell}&A_{ijk\ell}\end{pmatrix}.
\end{equation}
The product in equation~(\ref{prod}) is simplied
\begin{equation}
\sum_{i,k} \begin{pmatrix}K^P_iD_{ijk\ell}K^P_k & K^P_iC_{ijk\ell}\mathk_k\\ \mathk_i B_{ijk\ell}K^P_k&\mathk_iA_{ijk\ell}\mathk_k\end{pmatrix}    \begin{pmatrix}u^P_{\ell} \\ u^Q_{\ell}\end{pmatrix},
\end{equation}
which we rewrite as
\begin{equation}
\begin{pmatrix} {\bar D}_{j\ell} & {\bar C}_{j\ell}\\ {\bar B}_{j\ell}&{\bar A}_{j\ell}\end{pmatrix}    \begin{pmatrix}u^P_{\ell} \\ u^Q_{\ell}\end{pmatrix},
\end{equation}
where the bar indicates that sums over indices $i,k$ have been taken into account. 
Setting the high-frequency components to zero, we have three equations for each index $j$ and three unknowns for each index $\ell$ in 3-D space. Thus
\begin{equation}
\begin{pmatrix}\bara_{11} + \rho \omega^2 & \bara_{12} & \bara_{13} \\ \bara_{21} & \bara_{22} + \rho \omega^2 &\bara_{23}\\ \bara_{31}&\bara_{32}&\bara_{33}  +  \rho \omega^2\end{pmatrix}
\begin{pmatrix}u^Q_1\\u^Q_2\\u^Q_3\end{pmatrix} =
-\begin{pmatrix}\barb_{11} & \barb_{12} & \bara_{13} \\ \barb_{21} & \barb_{22} &\barb_{23}\\ \barb_{31}&\barb_{32}&\barb_{33}\end{pmatrix}
\begin{pmatrix}u^P_1\\u^P_2\\u^P_3\end{pmatrix} + \begin{pmatrix}f^Q_1\\f^Q_2\\f^Q_3\end{pmatrix},
\end{equation} 
and leaving aside the source for the moment, we obtain the corrector
\begin{equation}
-C A^{-1} B \begin{pmatrix}u^P_1\\u^P_2\\u^P_3\end{pmatrix},
\end{equation}
where
\begin{eqnarray}
A &=& \begin{pmatrix}\bara_{11} + \rho \omega^2 & \bara_{12} & \bara_{13} \\ \bara_{21} & \bara_{22} + \rho \omega^2 &\bara_{23}\\ \bara_{31}&\bara_{32}&\bara_{33}  +  \rho \omega^2\end{pmatrix},
\,\,\,\,\,B =\begin{pmatrix}\barb_{11} & \barb_{12} & \barb_{13} \\ \barb_{21} & \barb_{22} &\barb_{23}\\ \barb_{31}&\barb_{32}&\barb_{33}\end{pmatrix},\label{def.A}\\
\,\,\,\,C &=& \begin{pmatrix}\barc_{11} & \barc_{12} & \barc_{13} \\ \barc_{21} & \barc_{22} &\barc_{23}\\ \barc_{31}&\barc_{32}&\barc_{33}\end{pmatrix},\,\,\,\,\,
\,\,\,\,D = \begin{pmatrix}\bard_{11} & \bard_{12} & \bard_{13} \\ \bard_{21} & \bard_{22} &\bard_{23}\\ \bard_{31}&\bard_{32}&\bard_{33}\end{pmatrix}.
\end{eqnarray}
The effective wavespeed matrix is thus $ D - C A^{-1} B$, and is a function of temporal frequency because the definition of $A$ in equation~(\ref{def.A}) has an $\omega^2$ in it.
From this relation, each homogenized wavespeed component ${\bar\tau}_{ijkl}$ may be computing the sum above while keeping in mind that
$\bard_{j\ell} = \sum_{ik} K^P_i D_{ijk\ell}K^P_k$ and similarly for matrices $B$ and $C$. The difficulty lies
in being able to store and invert the matrix composed of $\bara_{ij}$. Computing matrix products may also be
non trivial. The homogenized source is
\begin{equation}
\begin{pmatrix} \bar{f}_1\\ \bar{f}_2\\ \bar{f}_3\end{pmatrix} = \begin{pmatrix}f^P_1\\f^P_2\\f^P_3\end{pmatrix}- C\,A^{-1}\begin{pmatrix}f^Q_1\\f^Q_2\\f^Q_3\end{pmatrix}.
\end{equation}

Taking a step back, we see that spatial homogenization creates a model of wavespeed that is a function of temporal frequency, which multiplies the (frequency-dependent) solution. This is a product in temporal Fourier domain, and therefore a convolution in the time domain. Thus squeezing the (convolutional) wavespeed operator in spatial Fourier domain results in broadening not just in real space but time as well. This is to be expected since eigenfunctions of the wave operator, which depend on both spatial wavenumber and temporal frequency form the basis on which wave solutions are projected. Reducing the size of this basis, i.e. discarding high-frequency solutions to the operator, has repercussions therefore in both the spatial and temporal domains, resulting in a broadened $(t,t')$, $(x,x')$. Ignoring the $(t,t')$ and $(x,x')$ dependencies introduce further errors when evaluating solutions. This analysis suggests the existence of an uncertainty principle that involves both space and time and may point to reasons why strongly interacting systems such as turbulence have not succumbed to renormalization group.

\end{document}